\documentclass[12pt]{amsart}
\usepackage{amssymb}
\usepackage[all,dvips,line,arrow,curve]{xy}

\newtheorem{thm}{Theorem}

\newtheorem{lem}{Lemma}
\newtheorem{cor}{Corollary}

\DeclareMathOperator{\Des}{Des}
\DeclareMathOperator{\des}{des}
\DeclareMathOperator{\fdes}{fdes}
\DeclareMathOperator{\ndes}{ndes}
\DeclareMathOperator{\maj}{maj}
\DeclareMathOperator{\fmaj}{fmaj}
\DeclareMathOperator{\nmaj}{nmaj}
\DeclareMathOperator{\wt}{wt}
\DeclareMathOperator{\st}{st}
\DeclareMathOperator{\eg}{neg}
\DeclareMathOperator{\flag}{flag}

\title{Euler-Mahonian distributions of type $B_n$}

\author[L. M. Lai]{Laurie M. Lai}

\author[T. K. Petersen]{T. Kyle Petersen}

\date{October 19, 2008}

\begin{document}

\maketitle

\begin{abstract}
Adin, Brenti, and Roichman introduced the pairs of statistics $(\ndes, \nmaj)$ and $(\fdes, \fmaj)$. They showed that these pairs are equidistributed over the hyperoctahedral group $B_n$, and can be considered ``Euler-Mahonian" in that they generalize the Carlitz identity. Further, they asked whether there exists a bijective proof of the equidistribution of their statistics. We give such a bijection, along with a new proof of the generalized Carlitz identity.
\end{abstract}

\section{Introduction}

Statistics such as inversion number, major index, and descent number have been well studied for the symmetric group $S_n$. MacMahon \cite{Mac} showed that inversion number and major index are equidistributed and Carlitz \cite{C} described the joint distribution of descent and major index---the ``Euler-Mahonian" distribution---with the following theorem.

\begin{thm}[\cite{C}]\label{thm:CarGes}
For positive integers $n$,
\begin{equation}\label{eq:CarGes}
\sum_{r\geq 0} [r+1]_q^n t^r = \frac{ \sum_{u \in S_n} t^{\des(u)} q^{\maj(u)} }{\prod_{i=0}^n (1-tq^i)},
\end{equation}
where $[r+1]_q = 1+q+\cdots+q^r$.
\end{thm}

The Coxeter group generalization of inversion number is \emph{length}. (An element $w$ in a Coxeter group has length $k$ if $w=s_1\cdots s_k$ is a minimal expression for $w$ as a product of simple reflections.) Adin and Roichman \cite{AR} generalized MacMahon's result to the hyperoctahedral group $B_n$ in demonstrating that a new statistic, the \emph{flag major} statistic, is equidistributed with length. Adin, Brenti, and Roichman \cite{ABR}, introduced new statistics on $B_n$, the \emph{negative descent}, the \emph{negative major}, and the \emph{flag descent}. They proved that the pairs of statistics ($\fdes$, $\fmaj$) and ($\ndes$, $\nmaj$) are equidistributed over $B_n$. Moreover, they showed these bivariate distributions are ``type $B_n$ Euler-Mahonian" in the sense that they generalize Theorem \ref{thm:CarGes} as follows. (Chow and Gessel provide an alternative generalization in \cite{ChowGes}.)

\begin{thm}\label{thm:ABR}
For positive integers $n$,
\begin{align}
\sum_{r\geq 0} [r+1]_q^n t^r &= \frac{ \sum_{w \in B_n} t^{\ndes(w)} q^{\nmaj(w)} }{(1-t)\prod_{i=1}^n (1-t^2q^{2i})} \quad \mbox{\cite[Theorem 3.2]{ABR}}, \label{eq:neg} \\
&= \frac{ \sum_{w \in B_n} t^{\fdes(w)} q^{\fmaj(w)} }{(1-t)\prod_{i=1}^n (1-t^2q^{2i})} \quad \mbox{\cite[Theorem 4.2]{ABR}.}\label{eq:flag}
\end{align}
\end{thm}

One of the drawbacks of \cite{ABR} is that while their proof of \eqref{eq:neg} is quite brief and elementary, their proof of \eqref{eq:flag} is rather indirect and tedious. Further, they remark (see the discussion after \cite[Corollary 4.5]{ABR}) that ``it would be interesting to have a direct combinatorial (i.e., bijective) proof" that these type $B_n$ Euler-Mahonian statistics are equidistributed. The purpose of this note is to provide such a proof (Theorem \ref{thm:bij}). Along the way we give new combinatorial proofs of \eqref{eq:neg} and \eqref{eq:flag} (Theorems \ref{thm:neg} and \ref{thm:flag}, respectively).

\section{Overview}

To any word $w=w_1\cdots w_n$ whose letters are totally ordered, let $\Des(w):=\{ i: w_i > w_{i+1}\}$, and let $\des(w):=|\Des(w)|$, called the \emph{descent set} and \emph{descent number}, respectively. We define the \emph{major statistic} as \[ \maj(w) := \sum_{i \in \Des(w)} i.\] Let $S_n$ denote the set of all permutations of the set $[n]:=\{1,2,\ldots,n\}$. The \emph{refined Eulerian polynomial} is the generating function for the joint distribution of $(\des,\maj)$ over $S_n$: \[ S_n(t,q) := \sum_{u\in S_n} t^{\des(u)}q^{\maj(u)}.\]
This function is the numerator of the right-hand side of \eqref{eq:CarGes}. For $q=1$, we have $S_n(t) = \sum_{u \in S_n} t^{\des(u)}$, the classical Eulerian polynomial shifted by $t^{-1}$.

The hyperoctahedral group, $B_n$, has as its elements all \emph{signed permutations} of $[n]$. For our purposes, these are all words $w=w_1\cdots w_n$ on the alphabet \[ [n]\cup [\bar n] = \{ \bar 1, 1, \bar 2, 2, \ldots, \bar n, n\} \] such that $|w|=|w_1|\cdots |w_n|$ is a permutation in $S_n$. (We write $\bar i$ for $-i$.) The precise way in which the flag and negative statistics are defined for $w$ in $B_n$ depends on how we choose to put a total order on $[n] \cup [\bar n]$, though the distribution over $B_n$ is independent of this choice. It will be most convenient to have the lexicographic order:
\[ \bar 1 < \bar 2 < \cdots < \bar n < 1 < 2 < \cdots < n. \] We remark that \cite{AR} uses this order, whereas \cite{ABR} uses the usual integer ordering.

Our approach is to provide two very explicit combinatorial proofs that the pairs of statistics in question are distributed over $B_n$ as \[ S_n(t,q)\cdot\prod_{i=1}^n (1+tq^i).\] For $q=1$ this is the intriguing formula $(1+t)^nS_n(t)$. (The first author in fact has a different combinatorial proof that flag descents are distributed in this fashion; see \cite{Lai}.)

In fact, we will demonstrate a finer result. Namely, that for each pair of statistics we have a way of assigning $2^n$ signed permutations to each unsigned permutation $u \in S_n$ so that the resulting collection has distribution \[t^{\des(u)}q^{\maj(u)}\prod_{i=1}^n(1+tq^i).\] The signed permutations in this collection are thus identified with a pair $(u,J)$, $J \subset [n]$, with ``weight" $t^{\des(u) + |J|}q^{\maj(u)+\sum_{j \in J} j}$. This weight-preserving bijection can then be composed to obtain the desired bijection. If $w$ corresponds to $(u,J)$ with respect to $(\ndes,\nmaj)$ and $v$ corresponds to $(u,J)$ with respect to $(\fdes,\fmaj)$, then we identify $w$ and $v$ with one another:
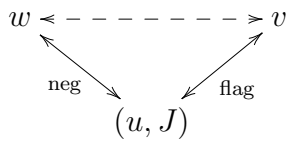
\[
      \xymatrix{
     w  \ar@{<-->}[rr]
       \ar@{<->}[dr]_{\eg} & & v           \\
       & (u,J) \ar@{<->}[ur]_{\flag} & }
\]
This idea is formalized with Theorem \ref{thm:bij}.

The paper is organized as follows. Section \ref{sec:neg} presents the negative statistics, Section \ref{sec:flag} discusses the flag statistics, and Section \ref{sec:bij} exhibits the desired bijection.

\section{The negative statistics}\label{sec:neg}

As in \cite{ABR}, define the \emph{negative descent number} of $w \in B_n$ as
\[ \ndes(w):= \des(w) + |\{i : w_i < 0\}|,\] and the
\emph{negative major index} as
\[ \nmaj(w):= \maj(w) + \sum_{w_i < 0} |w_i|.\]
For example if $w = \bar 5 7612\bar 4 \bar 3$, $\ndes(w) = 4 +3 = 7$, and $\nmaj(w) = 16+ 12 = 28$. Let \[B^{(\eg)}_n(t,q) := \sum_{w \in B_n} t^{\ndes(w)}q^{\nmaj(w)},\] be the generating function for this pair of statistics.

To any $u=u_1\cdots u_n$ in $S_n$ we can associate $2^n$ signed permutations in the following manner.
Define the \emph{standardization} of a signed permutation $w \in B_n$, $\st(w)$, to be the unsigned permutation in $S_n$ that is obtained by replacing the smallest letter of $w$ with $1$, the next smallest with $2$, and so on. For example, $\st(\bar 5 7612\bar 4 \bar 3) = 3764521$.

Now given any $u$ in $S_n$, let \[B(u) = \{ w \in B_n : \st(w) = u\}.\] We have \[ B_n = \bigcup_{u \in S_n} B(u) \quad \mbox{ (disjoint union). }\]

\medskip

\noindent\textbf{Remark.}
This partition of $B_n$ is employed by Adin, Brenti, and Roichman (see equation (5) of \cite{ABR}), though in different language. It plays a role in their proof of their Theorem 3.2 (equation \eqref{eq:neg} above) and a refined version of the same result from a different work \cite[Theorem 6.7]{ABR2}. This partition will play a role in our proof as well, though the proofs are fundamentally distinct.

\medskip

It is easy to see that every element of $B(u)$ is uniquely determined by its set of negative letters. Let $u_J$ denote the member of $B(u)$ that is a permutation of the set \[ \{ i : i \in [n]\setminus J\} \cup \{ \bar j : j \in J\}.\] For example, if $u = 3142$, then $u_{\{1,3\}} = 2\bar 1 4 \bar 3$ since $\bar 1$ and $\bar 3$ are the negative letters and $\st(2\bar 1 4 \bar 3) = 3142$.

The following result and Theorem \ref{thm:CarGes} implies equation \eqref{eq:neg} of Theorem \ref{thm:ABR}.

\begin{thm}\label{thm:neg}
We have \[ B_n^{(\eg)}(t,q) = S_n(t,q)\cdot\prod_{i=1}^n (1+tq^i).\]
\end{thm}

Theorem \ref{thm:neg} follows immediately from the subsequent stronger claim.

\begin{thm}\label{thm:refineneg}
For any $u \in S_n$, we have
\[ t^{\des(u)}q^{\maj(u)}\prod_{i=1}^n(1+tq^i) = \sum_{w \in B(u)} t^{\ndes(w)}q^{\nmaj(w)}.\]
\end{thm}

\begin{proof}
Let $w = u_J$ for any $J \subset [n]$. We want to show that \[ t^{\ndes(w)}q^{\nmaj(w)} = t^{\des(u)}t^{\maj(u)}\prod_{j \in J} tq^j.\]

Since $\st(w) = u$, we have $\Des(w) = \Des(u)$. In particular, $\des(w) = \des(u)$ and $\maj(w) = \maj(u)$. By definition of $w$, we know $\{ w_i: w_i < 0\} = J$. Comparing with the definition of the negative statistics, this completes the proof.
\end{proof}

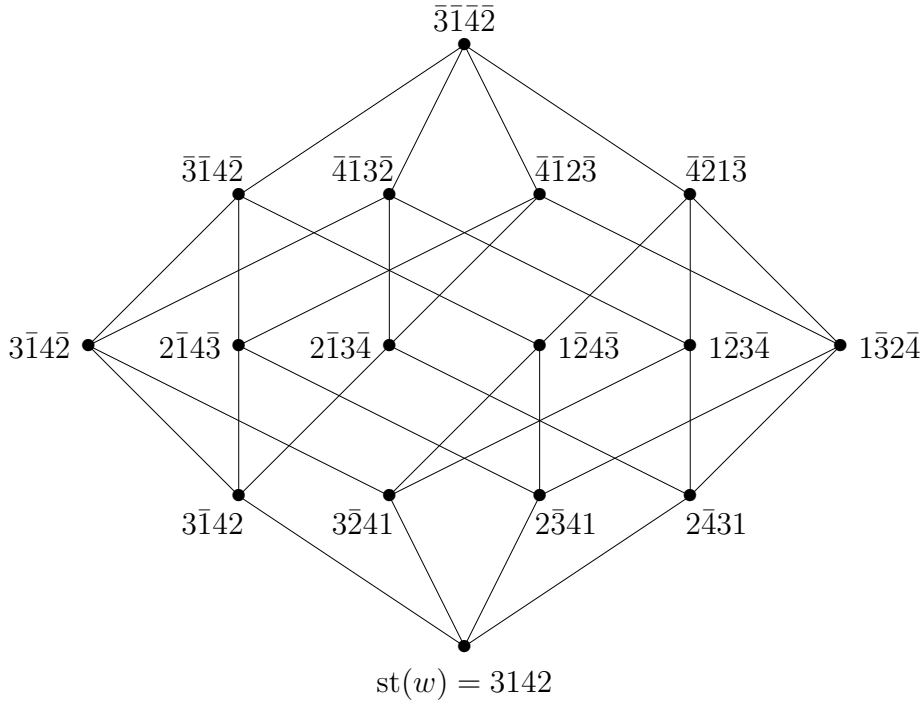
\begin{figure}[h]
\begin{xy}
0;<1cm,0cm>:
(-3,0)*{}; (5,0); (2,2) **@{-}, (4,2) **@{-}, (6,2) **@{-}, (8,2) **@{-}, (2,2); (0,4) **@{-}, (2,4) **@{-},(4,4) **@{-}, (4,2); (0,4) **@{-}, (6,4) **@{-}, (8,4) **@{-}, (6,2); (2,4) **@{-}, (6,4) **@{-}, (10,4) **@{-}, (8,2); (4,4) **@{-}, (8,4) **@{-},(10,4) **@{-}, (2,6); (0,4) **@{-}, (2,4) **@{-}, (6,4) **@{-}, (4,6); (0,4) **@{-}, (4,4) **@{-},  (8,4) **@{-}, (6,6); (2,4) **@{-}, (4,4) **@{-},  (10,4) **@{-}, (8,6); (6,4) **@{-}, (8,4) **@{-},  (10,4) **@{-}, (5,8); (2,6) **@{-}, (4,6) **@{-},  (6,6) **@{-}, (8,6) **@{-},
(5,0)*{\bullet}, (2,2)*{\bullet}, (4,2)*{\bullet}, (6,2)*{\bullet}, (8,2)*{\bullet}, (0,4)*{\bullet}, (2,4)*{\bullet}, (4,4)*{\bullet}, (6,4)*{\bullet}, (8,4)*{\bullet}, (10,4)*{\bullet}, (2,6)*{\bullet}, (4,6)*{\bullet}, (6,6)*{\bullet}, (8,6)*{\bullet}, (5,8)*{\bullet},  (5,-.5)*{\st(w) = 3142}, (1.65,1.65)*{ 3 \bar 1 4 2}, (3.65,1.65)*{ 3 \bar 241}, (6.35,1.65)*{2 \bar 3 41}, (8.35,1.65)*{2 \bar 4 3 1}, (-.65,4)*{3 \bar 1 4 \bar 2}, (1.35,4)*{2 \bar 14 \bar 3}, (3.35,4)*{2 \bar 1 3 \bar4}, (6.65,4)*{1 \bar 2 4 \bar 3}, (8.65,4)*{1 \bar 2 3 \bar 4}, (10.65,4)*{1 \bar 3 2 \bar 4}, (1.65,6.35)*{\bar 3\bar 1 4 \bar 2}, (3.65,6.35)*{\bar 4\bar 1 3 \bar2}, (6.35, 6.35)*{\bar 4\bar12\bar3}, (8.35, 6.35)*{\bar 4 \bar 21 \bar3}, (5, 8.35)*{\bar 3 \bar 1\bar4\bar2}
\end{xy}
\caption{$B(u)$ for $u = 3142$}\label{fig:ex2}
\end{figure}

It is informative to observe that Theorem \ref{thm:refineneg} implies the set $B(u)$ has, up to a shift, a weight-preserving bijection with a boolean algebra given by \[ J \longleftrightarrow u_J. \] 
That is, if $w \in B(u)$ corresponds to $J = \{ j_1, \ldots, j_k\}$, then $\ndes(w) = \des(u)+k$ and $\nmaj(w) = \maj(u) + (j_1 + \cdots +j_k)$.
Figure \ref{fig:ex2} provides an illustration for $u =3142$.

\section{The flag statistics}\label{sec:flag}

From \cite{ABR}, our definitions for the flag statistics are as follows. For $w \in B_n$, \[ \fdes(w) := \begin{cases}
 2\des(w) + 1 & \mbox{if } w_1 < 0 \\
 2\des(w) & \mbox{if } w_1 > 0,
 \end{cases}
 \]
and \[ \fmaj(w) := 2\maj(w) + |\{ i : w_i < 0\}|.\]
So, for example, if $w = \bar 5 7612\bar 4 \bar 3$, then $\Des(w) = \{ 2, 3, 5, 6\}$, $\fdes(w) = 2\cdot 4 + 1 = 9$, and $\fmaj(w)= 2\cdot 16 + 3 = 35$. Let $\wt(w) = t^{\fdes(w)}q^{\fmaj(w)}$ be the \emph{weight} of $w$. The generating function is: \[ B_n^{(\flag)}(t,q) = \sum_{w \in B_n} \wt(w) = \sum_{w \in B_n} t^{\fdes(w)}q^{\fmaj(w)}.\]

Before proceeding, it will be helpful to provide another partition of $B_n$. To any $u = u_1\cdots u_n$ in $S_n$, we can assign $2^n$ signed permutations by independently assigning minus signs to the letters of $u$ in all possible ways. That is, let \[ B'(u) = \{ w \in B_n : |w| = u \}.\] We have \[ B_n = \bigcup_{u \in S_n} B'(u) \quad \mbox{ (disjoint union).}\]

Just as with Theorem \ref{thm:neg}, we have the following, which, with Theorem \ref{thm:CarGes} completes the proof of Theorem \ref{thm:ABR}.

\begin{thm}[\cite{ABR}, Theorem 4.4]\label{thm:flag}
We have \[ B_n^{(\flag)}(t,q) = S_n(t,q)\cdot\prod_{i=1}^n (1+tq^i) .\]
\end{thm}

\begin{cor}[\cite{ABR}, Corollary 4.5]\label{cor:equidist}
The statistics $(\fdes,\fmaj)$ and $(\ndes,\nmaj)$ are equidistributed, i.e., \[ B^{(\flag)}_n(t,q) = B_n^{(\eg)}(t,q).\]
\end{cor}

Theorem \ref{thm:flag} will follow from the following theorem, analogous to Theorem \ref{thm:neg}. With the combinatorial proof of Theorem \ref{thm:refined} in hand we will be able to provide the desired bijective proof of Corollary \ref{cor:equidist}.

\begin{thm}\label{thm:refined}
For any $u$ in $S_n$, we have:
\[ t^{\des(u)}q^{\maj(u)}\prod_{i=1}^n (1+tq^i) = \sum_{w \in B'(u)} t^{\fdes(w)}q^{\fmaj(w)}.\]
\end{thm}

Theorem \ref{thm:refined} will follow from Lemma \ref{lem:symdiff}. First, we need one more idea. Let $\Delta_i$ be the operator that negates the first $i$ letters of a signed permutation, i.e.,
\[ \Delta_i w = \overline{w_1\cdots w_i} w_{i+1}\cdots w_n.\]

\begin{lem}\label{lem:flip}
Fix $u \in S_n$ and let $w \in B'(u)$ be such that all the letters $w_1,\ldots, w_{j+1}$ have the same sign. Then \[ \wt(\Delta_j w) =
\begin{cases}
\wt(w)/tq^j & \mbox{if } j \in \Des(u)\\
\wt(w)\cdot tq^j & \mbox{if } j \notin \Des(u).
\end{cases}
\]
\end{lem}

\begin{proof}
Suppose everything to the left of $w_{j+1}$ (inclusive) is positive. Then if $j \notin \Des(u)$ we have lost nothing and gained $tq^j$, since $w_1$ is now negative and we have $j$ new negative numbers. If we have $j \in \Des(u)$, then we still gain $tq^j$, but we have lost the descent in $j$, and with it $t^2q^{2j}$, for a net weight change of $t^{-1}q^{-j}$.

If everything to the left of $w_{j+1}$ is negative the situation is similar. If $j \in \Des(u)$ we have gained no new descents, lost $j$ negative numbers, and changed $w_1$ from negative to positive; a net weight change of $t^{-1}q^{-j}$. If we have $j \notin \Des(u)$, then we gain $t^2q^{2j}$ as $j$ becomes a descent, but we have lost $j$ negative numbers and $w_1$ has gone from negative to positive, yielding a change of $t^2q^{2j}/tq^j = tq^j$.

This completes the proof.

\end{proof}

Notice that the $\Delta_i$ commute, so for a subset $J = \{ j_1 < \cdots < j_k\}$ of $[n]$, there is no ambiguity in defining the composite operator $\Delta_J = \Delta_{j_1} \cdots \Delta_{j_k}$.

\begin{lem}\label{lem:symdiff}
For any $u \in S_n$, we have
\[ \wt(\Delta_J u) = t^{\des(u)}q^{\maj(u)}\prod_{j \in J \vartriangle \Des(u)} tq^j,\] where $J \vartriangle \Des(u) = (J \cup \Des(u)) \setminus (J \cap \Des(u))$ is the symmetric difference of $J$ and $\Des(u)$.
\end{lem}

\begin{proof}
First notice that $\fdes(u) = 2\des(u)$ and $\fmaj(u) = 2\maj(u)$. Since the $\Delta_j$ commute, we can apply $\Delta_{j_k}$ first, $\Delta_{j_{k-1}}$ second, and so on. This allows us to apply Lemma \ref{lem:flip} repeatedly, giving us
\begin{align*}
 \wt(\Delta_J u) &= t^{2\des(u)}q^{2\maj(u)} \frac{\prod_{j \in J\setminus \Des(u)} tq^j}{\prod_{i \in J\cap \Des(u)} tq^i}\\
&= t^{\des(u)}q^{\maj(u)} \prod_{j \in J \vartriangle \Des(u)} tq^j,
\end{align*}
as desired.
\end{proof}

Now we can prove Theorem \ref{thm:refined}.

\begin{proof}[Proof of Theorem \ref{thm:refined}]
By Lemma \ref{lem:symdiff}, it suffices to show that \[ B'(u) = \{ \Delta_J u : J \subset [n]\}.\] Because the two sets have the same cardinality, it suffices to show $\Delta_J u \in B'(u)$ for any $J$. This is clear since $|\Delta_J u| = u$.
\end{proof}

\section{The bijection}\label{sec:bij}

From Lemma \ref{lem:symdiff}, we notice in particular that there is a unique signed permutation with the same flag statistics as $u$ has ordinary statistics, obtained by taking $J = \Des(u)$. Define \[ w_{\min}(u) = \Delta_{\Des(u)}u.\]

Given any signed permutation $w$ we have $\Delta_i^2 w = w$ and, if $i < j$, \[ \Delta_i \Delta_j w = \Delta_j\Delta_i w = w_1 \cdots w_i \overline{w_{i+1} \cdots  w_j} w_{j+1} \cdots w_n.\] In particular the composite operator $\Delta_{i}\Delta_{i+1}$ simply negates letter $i+1$. Since the  $\Delta_i$ commute and generate all possible sign changes, we can identify permutations in $B'(u)$ with subsets of $[n]$ by the correspondence
\[ J \longleftrightarrow \Delta_J w_{\min}(u) = \Delta_{J\vartriangle \Des(u)} u. \]

Lemma \ref{lem:symdiff} can now be interpreted to  show that the boolean algebra generated by subsets of $[n]$ respects the statistics $(\fdes,\fmaj)$ on $B'(u)$. Specifically, if  $w \in B'(u)$ corresponds to $J = \{ j_1, \ldots, j_k\}$, then $\fdes(w) = \des(u)+k$ and $\fmaj(w) = \maj(u) + (j_1 + \cdots +j_k)$.

\begin{figure}[h]
\begin{xy}
0;<1cm,0cm>:
(-3,0)*{}; (5,0); (2,2) **@{-}, (4,2) **@{-}, (6,2) **@{-}, (8,2) **@{-}, (2,2); (0,4) **@{-}, (2,4) **@{-},(4,4) **@{-}, (4,2); (0,4) **@{-}, (6,4) **@{-}, (8,4) **@{-}, (6,2); (2,4) **@{-}, (6,4) **@{-}, (10,4) **@{-}, (8,2); (4,4) **@{-}, (8,4) **@{-},(10,4) **@{-}, (2,6); (0,4) **@{-}, (2,4) **@{-}, (6,4) **@{-}, (4,6); (0,4) **@{-}, (4,4) **@{-},  (8,4) **@{-}, (6,6); (2,4) **@{-}, (4,4) **@{-},  (10,4) **@{-}, (8,6); (6,4) **@{-}, (8,4) **@{-},  (10,4) **@{-}, (5,8); (2,6) **@{-}, (4,6) **@{-},  (6,6) **@{-}, (8,6) **@{-},
(5,0)*{\bullet}, (2,2)*{\bullet}, (4,2)*{\bullet}, (6,2)*{\bullet}, (8,2)*{\bullet}, (0,4)*{\bullet}, (2,4)*{\bullet}, (4,4)*{\bullet}, (6,4)*{\bullet}, (8,4)*{\bullet}, (10,4)*{\bullet}, (2,6)*{\bullet}, (4,6)*{\bullet}, (6,6)*{\bullet}, (8,6)*{\bullet}, (5,8)*{\bullet},  (5,-.5)*{w_{\min} = 3\bar 1\bar4 2}, (1.65,1.65)*{\bar 3 \bar 1 \bar 4 2}, (3.65,1.65)*{\bar 3 1 \bar 4 2}, (6.35,1.65)*{\bar 3 1 4 2}, (8.35,1.65)*{\bar 3 14 \bar 2}, (-.65,4)*{31\bar 4 2}, (1.35,4)*{3142}, (3.35,4)*{314\bar 2}, (6.65,4)*{3\bar 1 42}, (8.65,4)*{3\bar 1 4\bar 2}, (10.65,4)*{3\bar 1 \bar 4\bar 2}, (1.65,6.35)*{\bar 3\bar 1 42}, (3.65,6.35)*{\bar 3\bar 1 4 \bar2}, (6.35, 6.35)*{\bar 3\bar1\bar4\bar2}, (8.35, 6.35)*{\bar 3 1\bar4\bar2}, (5, 8.35)*{31\bar4\bar2}
\end{xy}
\caption{$B'(u)$ for $u = 3142$}\label{fig:ex}
\end{figure}
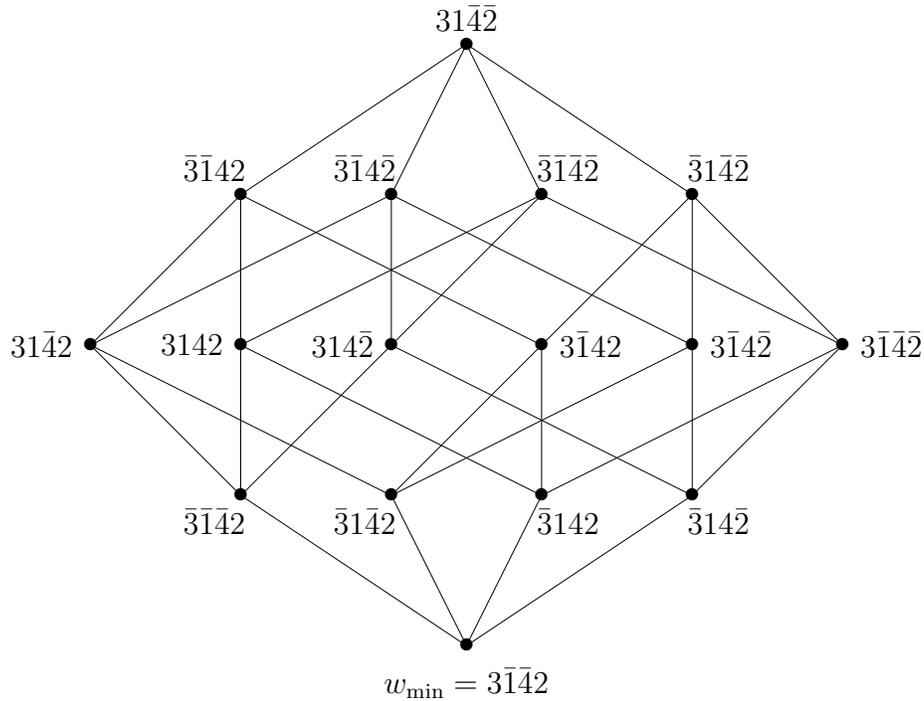

See Figure \ref{fig:ex} for an example with $u=3142$.

Implicit in here is the bijection that Adin, Brenti, and Roichman wanted.

\begin{thm}\label{thm:bij}
The bijection \[ u_J \longleftrightarrow  \Delta_{J\vartriangle \Des(u)} u\] is weight-preserving in that
\begin{align*}
\ndes(u_J) &= \fdes(\Delta_{J\vartriangle \Des(u)} u) \\
\nmaj(u_J) &= \fmaj(\Delta_{J\vartriangle \Des(u)} u).
\end{align*}
\end{thm}

\noindent\textbf{Example.} Let $w = \bar 2 1 \bar 3 \bar 4$ with $\ndes(w)=4$ and $\nmaj(w)=11$. We first find $u=\st(w)=1423$. Since $J=\{2,3,4\}$ is the set of negative numbers in $w$ and $\Des(u)=\{2\}$, we take $v=\Delta_3 \Delta_4 ( 1 4 2 3) =1 4 2 \bar 3$, with the desired flag descent and flag major numbers, $\fdes(v) = 4$ and $\fmaj(v)=11$.

In the other direction, let $v=43\bar1\bar2$ with $\fdes(v)=4$ and $\fmaj(v)=8$. First we see that $u=|v| = 4312$ and $v = \Delta_2\Delta_4 u$. Since $\Des(u) = \{1,2\}$, we get $J = \{1,4\}$ and $u_J = 32\bar 1 \bar 4$ with $\ndes(u_J) = 4$ and $\nmaj(u_J) = 8$, as desired.

\end{document}